\documentclass[10pt,reqno,letterpaper,oneside]{amsart}
\usepackage{amssymb}
\usepackage[USenglish]{babel}
\usepackage{dsfont} 
\usepackage{mathrsfs}   
\usepackage{paralist} 
\usepackage[nice]{nicefrac}   
\usepackage{ifthen}
\newif\ifpdf
\ifx\pdfoutput\undefined
   \pdffalse
\else
   \pdfoutput=1
   \pdftrue
\fi

\ifpdf
\usepackage{microtype}
\usepackage[pdftex]{color} 
\usepackage[pdftex]{geometry}
\usepackage[pdftex]{hyperref}
\else %
\usepackage{color}
\usepackage[hypertex]{hyperref}
\fi

\definecolor{rubrojo}{rgb}{0.82,0.00,0.41}
\definecolor{rubceles}{rgb}{0.00,0.67,0.67}
\definecolor{rubverde}{rgb}{0.00,0.62,0.00}
\definecolor{rubblue}{rgb}{0.00,0.00,0.70}

\ifpdf
\hypersetup{%
       pdftitle={MHD equations without diffusion},
       pdfauthor={Agapito Schonbek},
       pdfstartview={FitH},plainpages=false,
       bookmarksnumbered,breaklinks=true,
       urlcolor=rubblue,
       menucolor=rubblue,
       citecolor=rubblue,
       linkcolor=rubblue,
       colorlinks=true}
\else
\hypersetup{%
       urlcolor=rubblue,
       menucolor=rubceles,breaklinks=true,
       citecolor=rubblue,
       linkcolor=rubblue,
       colorlinks=true}
\fi

\definecolor{deficolor}{rgb}{0.85,0.97,0.99}
\definecolor{teocolor}{rgb}{0.80,1.00,0.80}

\theoremstyle{plain} 
\newtheorem{teo}{Theorem}
\newtheorem{cor}[teo]{Corollary}
\newtheorem{prop}{Proposition}

\newtheorem*{teononum}{Theorem}

\theoremstyle{definition} %

\theoremstyle{remark} %

\theoremstyle{plain} %

\numberwithin{equation}{section}

\DeclareMathOperator{\diver}{div}%

\newcommand{\R}{\mathds{R}}        

\newcommand{\eps}{\varepsilon}     

\newcommand{\norm}[1]{\left\lVert#1\right\rVert}  
\newcommand{\abs}[1]{\left\lvert#1\right\rvert}   

\newcommand{\beqst}{\begin{equation*}}
\newcommand{\eeqst}{\end{equation*}}
\newcommand{\beq}{\begin{equation}}
\newcommand{\eeq}{\end{equation}}

\begin{document}
\title[Non-uniform Decay of MHD Equations]{Non-uniform Decay of MHD equations \\%
with and without magnetic diffusion}
\author{Rub\'en~Agapito}
\address{Dept. of Mathematics, UCSC, CA 95064}
\email{ragapito@math.ucsc.edu}
\urladdr{http://math.ucsc.edu/~ragapito/}

\author{Maria Schonbek$^*$}
\email{schonbek@math.ucsc.edu}
\urladdr{http://math.ucsc.edu/~schonbek/}
\thanks{The work of M. Schonbek was partially supported
  by NSF grant DNS-0600692}

\begin{abstract}
We consider the long time behavior of solutions to the
magnetohydrodynamics equations in two and three spatial dimensions.
It is shown that in the absence of magnetic diffusion, if strong
bounded solutions were to exist their energy cannot present any
asymptotic oscillatory behavior, the diffusivity of the velocity is
enough to prevent such  oscillations. When magnetic diffusion is
present  and the data  is only in $L^2$, it is shown that the
solutions decay to zero without a rate, and this nonuniform decay is
optimal.
\end{abstract}

\maketitle %
\section{Introduction}

We consider the Magnetohydrodynamics equations (MHD) in two and
three dimensions. We deal with questions regarding long time
behavior of solutions to the MHD with and without magnetic
diffusion. The MHD equations model the interactions between a
magnetic field and a viscous incompressible fluid of moving
electrically charged particles.

In non-dimensional form the equations can be expressed by
\begin{equation}\label{primero}
\begin{split}
\dfrac{\partial}{\partial t}u +(u\cdot\nabla)u-S(B\cdot\nabla)B
+\nabla (P+\frac{S}{2}\abs{B}^2) & =\frac{1}{\mathrm{Re}}\Delta u, \\
\dfrac{\partial}{\partial t}B +(u\cdot\nabla)B -(B\cdot\nabla)u
&=\delta \Delta B,\\
\nabla\cdot u &=0,\,\nabla\cdot B=0,\\
 u(x,0) &=u_0(x),\,B(x,0)=B_0(x),
\end{split}
\end{equation}
where $u=u(x,t)=(u^1(x,t),u^2(x,t),\dotsc,u^n(x,t)), B=B(x,t)$ and
$P=P(x,t)$ denote the unknown velocity, the magnetic field and
pressure of the fluid at a point $(x,t)\in\R^n\times\R_+$,
respectively. The term $\frac{\abs{B}^2}{2}$ denotes magnetic
pressure. The positive constants appearing in the equations are
$\mathrm{Re}$, the Reynolds number, $\mathrm{Rm}$, the magnetic
Reynolds number, and $S=M^2/(\mathrm{Re}\,\mathrm{Rm})$, where $M$
is the Hartman number. For the sake of notational simplicity, and
with minor loss of generality, we set all these constants equal to
one. After rescaling $u$ and $B$, let $p=P+\frac{1}{2}S\abs{B}^2$
denote the total pressure, Equation (\ref{primero}) can be rewritten
as
\begin{align}
\dfrac{\partial}{\partial t}u +(u\cdot\nabla)u-(B\cdot\nabla)B
+\nabla p & =\Delta u, \label{mhd1}\\
\dfrac{\partial}{\partial t}B +(u\cdot\nabla)B -(B\cdot\nabla)u
&=\delta\Delta B, \label{mhd2}\\
\nabla\cdot u &=0,\,\nabla\cdot B=0, \label{mhd3}\\
 u(x,0) &=u_0(x),\,B(x,0)=B_0(x).\label{mhd4}
\end{align}
The initial data $(u_0(x),B_0(x))$ will be chosen below in
appropriate spaces. Derivations of these equations can be found in
\cite{cha,cow,landau}.

Many authors have studied MHD equations from the point of view of
existence and long time behavior. Without making a complete list of
all authors we would like to mention some of the relevant
literature. In particular important results on existence were
obtained, among others, in \cite{kozo87,xin05,xin052,wu02}. In the
direction of decay interesting results can be found in the papers
\cite{india89,kim02,maria96}. The methods used for decay in
\cite{india89} were based on Fourier Splitting \cite{maria85}. The
paper \cite{kim02} uses Gevrey regularity and it is based on ideas
developed by Foias and Temam. Similar methods were used for the
Navier-Stokes equations and can be found in \cite{oliver00}.

Several questions will be addressed regarding the long time behavior
of the solutions. In this regard, in the absence of magnetic
diffusion, we are going to analyze the following problem.

\emph{Is the diffusion introduced by the velocity alone sufficient
to prevent compensatory oscillations?} Specifically, simple
calculations shows that the following energy inequality holds when
$\delta=0$,
\beqst
\dfrac{d}{dt}\bigl(\norm{u(t)}^2_2 +\norm{B(t)}^2_2\bigr) \le -2
\norm{\nabla u(s)}^2_2.
\eeqst
This shows that the combined energy decays, but allows the
possibility of separate oscillations in $\norm{u}_2$ and
$\norm{B}_2$ that could compensate each other. In this paper we show
that in the case that there were to exist strong solutions which are
bounded such oscillations can never occur. Specifically it is shown
that

\begin{teo}\label{teore1}
Let $n=3$. Suppose that there exists $(u,B)$ a strong bounded
solution to the MHD equations with $\delta=0$, magnetic field
bounded in $L^{\infty}$, and data $(u_0,B_0)\in (L^1\cap L^2\times
L^2\cap L^{\infty})$. Then
\beqst
\norm{u(t)}_2\to 0,\quad\text{and}\quad \norm{B(t)}_2\to M,
\eeqst
as $t\to\infty$, where $M$ is some positive constant.
\end{teo}

The previous theorem shows that the diffusion in the velocity is
sufficient to prevent compensate oscillations between the two
energies.

The next question we addressed is in regards to decay of solutions
to MHD equations with diffusion both in the velocity and the
magnetic field. Here it is supposed that the data is only in $L^2$
space. In this case it is shown that the energies of the velocity
and the magnetic field decay to zero without a rate. Moreover, it is
shown that this is optimal, that is, that cannot be a uniform rate
for the energy of the solutions with data exclusively in $L^2$. We
show that

\begin{teo}\label{teore3}
Let $n=2,3$. If $(u,B)$ is a weak solution to the MHD equations with
$\delta>0$ and data $(u_0,B_0)\in \bigl(L^2(\R^n)\bigr)^2$, then
\beqst
\lim_{t\to\infty} \left(\norm{u(t)}^2_2 + \norm{B(t)}^2_2\right)=0.
\eeqst
\end{teo}

The proof of this theorem is given first formally. To make the
result rigorous we apply the method to smooth approximations and
then pass to the limit.

With regard to the optimality of this last result we show that

\begin{teo}
There exist no functions $G(t,\beta)$ and $H(t,\gamma)$ with the
following two properties. If $(u,B)$ is a solution to equations
(\ref{mhd1})-(\ref{mhd4}) with $\delta>0$ and data $(u_0,B_0)\in
(L^2(\R^n))^2$, $n=2,3$, then
\begin{enumerate}[i{)}]
\item $\norm{u(t)}_2\le G(t,\norm{u_0}_2)$,
\item $\lim_{t\to 0}G(t,\beta)=0$, for all $\beta >0$.
\end{enumerate}
\end{teo}

\vspace{1ex}%
The last part of the paper focuses on extending Kato's pioneering
work on $L^p$ decay for Navier-Stokes equations \cite{kato84} to the
MHD equations with magnetic diffusion. In particular we note that a
simple modification of Kato's work \cite{kato84} yields equivalent
results for the MHD equations, then combined with our result on
non-uniform decay in $L^2$ gives a slight improvement of the decay
rates.

\begin{cor}
There is $\lambda >0$ such that for $\norm{u_0}_2\le \lambda$ the
global solution of the equation (\ref{ecuateo4}) for $q\ge m,$ and
for $2\le r\le q$
\beqst
\lim_{t\to\infty} t^{\frac{r-2}{2r}}\norm{(u,B)}_r=0
\eeqst
\end{cor}

\section{Notation}

We denote
\begin{align*}
C^{\infty}_{0,\sigma}=C^{\infty}_{0,\sigma}(\R^n)\; : & \text{ space
  of infinitely differentiable functions with compact}\\ &
  \text{ support and divergence free,}\\
L^2_{\sigma}=L^2_{\sigma}(\R^n)\; :& \text{ completion of
  }C^{\infty}_{0,\sigma}\text{ in the }L^2\text{-norm }\norm{\cdot}_2,\\
\dot{H}^1_{\sigma}=\dot{H}^1_{\sigma}(\R^n)\; : & \text{ completion
  of }C^{\infty}_{0,\sigma}\text{ in the homogeneous }H^1\text{-norm
  }\norm{\nabla\cdot}_2.
\end{align*}
The Fourier transform of $\varphi$ will be denoted by
$\mathscr{F}\{\varphi \}=\Hat{\varphi}=\int_{\R^n} e^{-i\xi
x}\varphi(x)\,dx$ and its inverse transform by
$\mathscr{F}^{-1}\{\varphi \}=\check{\varphi}=\frac{1}{2\pi^n}\int
e^{ix\xi}\Hat{\varphi}(\xi)d\xi$. Also,
\beqst
L^p(a,b;L^q)=\left\{ f:(a,b)\times\R^n\rightarrow\R^n:
\norm{f}_{L^p(a,b;L^q)}=\Biggl( \int^b_a \norm{f(\tau)}^p_q~d\tau
\Biggr)^{\nicefrac{1}{p}}<\infty \right\}.
\eeqst
The notation $\norm{\cdot}_{L^{p,q}}$ will be used for the norm of
$L^p(0,\infty;L^q)$, and $\langle f,g \rangle=\int f g\,dx$ for the
inner product in $L^2$. Various constants are simply denoted by $C$.

\section{Preliminary Results}

In this section auxiliary results which will be needed in the sequel
are obtained. We consider the MHD equations with no magnetic
diffusion ($\delta=0$). Some of the results presented are standard
and their proofs are included in the appendix for completeness. The
results below follow ideas of \cite{ors97}.

We start with an estimate for the Fourier transform of the velocity.

\begin{prop}\label{decay}
 Let $(u,B)$ be a mild solution to the MHD equations (\ref{mhd1})-(\ref{mhd4}).
 Assume the initial data $u_0, B_0$ is in $L^1(\R^3)\cap L^2(\R^3)$.
 Then
 \beqst
\abs{\Hat{u}(t)}\le C\biggl(1+\dfrac{1}{\abs{\xi}} \biggr),
 \eeqst
 where $C$ is a constant.
\end{prop}

\noindent{\bf Proof.} See Appendix, Proposition \ref{proofdecay}.
\hfill$\square$

The proofs in this section are formal. To make them rigorous is
suffices to apply them to approximations using retarded
mollifications such as the ones constructed by \cite{ckn82,leray34},
and \cite{xin05} for the MHD equations and then pass to the limit.
For details regarding our proofs see \cite{ors97} were the procedure
has been done for the solutions to the Navier-Stokes equations. The
extension to MHD is straightforward.

We first recall a standard energy inequality
\beqst
\norm{u(t)}^2_2 +\norm{B(t)}^2_2 +2\int^t_0 \norm{\nabla
u(s)}^2_2~ds =\norm{u_0}^2_2 +\norm{B_0}^2_2.
\eeqst
This follows easily by multiplying (\ref{mhd1}) by $u$, (\ref{mhd2})
by $B$, adding the equations, and then integrating in space and
time.

The next proposition gives a generalized energy inequality.

\begin{prop}
Assume $(u,B)$ is a weak solution of Equations
(\ref{mhd1})-(\ref{mhd4}). For $E(t)\in C^1(\R;\R_+)$ with $E(t)\ge
0$ and $\psi\in C^1(\R;C^1\cap L^2)$, the weak solution satisfies
\beq\label{genener}
\begin{split}
E(t)\norm{\psi*u(t)}^2_2 & = E(s)\norm{\psi*u(s)}^2_2 + \int^t_s
  E'(\tau)\norm{\psi*u(\tau)}^2_2 ~d\tau \\
 & +2\int_s^t E(\tau)\Bigl(\langle \psi' *u(\tau),\psi*u(\tau)\rangle
   -\norm{\nabla\psi*u(\tau)}^2_2\Bigr) ~d\tau \\
 & -2\int^t_s E(\tau)\bigl(\langle u\cdot\nabla u(\tau),
    \psi*\psi*u(\tau)\rangle \bigr. \\
 &  \phantom{-2\int^t_s E(\tau)}+\bigl. \langle B\cdot\nabla B(\tau),
    \psi*\psi*u(\tau)\rangle \bigr) ~d\tau
\end{split}
\eeq
\end{prop}

\noindent{\bf Proof.} Multiply Equation (\ref{mhd1}) by
$E(t)\psi*\psi*u(t)$ and integrate by parts to get
\begin{align*}
\tfrac{d}{dt}\bigl(E(t)\norm{\psi*u(t)}^2_2\bigr) & =
E'(t)\norm{\psi*u(t)}^2 +2
    E(t)\bigl\{\langle \psi'*u(t),\psi*u(t)\rangle-\norm{\nabla\psi
    *u(t)}^2_2 \bigr\} \\
    & \phantom{=}~-2E(t)\bigl\{ \langle u\cdot\nabla u,\psi*\psi*u(t)\rangle
      -\langle B\cdot\nabla B,\psi*\psi*u(t)\rangle\bigr\}
\end{align*}
Integrating the preceding equation in the interval $(s,t)$ yields
Equation (\ref{genener}).\hfill$\square$

Corollaries \ref{corgen} and \ref{corgen2} follow as an easy
consequence.

\begin{cor}\label{corgen}
Let $(u,B)$ be a weak solution of (\ref{mhd1})-(\ref{mhd4}). Let
$\varphi\in L^2(\R^3)$, then
\begin{multline}\label{mylab}
\norm{\check{\varphi}* u(t)}_2^2 \le
\norm{e^{\Delta(t-s)}\check{\varphi}* u(s)}_2^2 + 2\int^t_s\Bigl(
 \bigl|\langle u\cdot\nabla u, e^{2\Delta(t-\tau)}\check{\varphi}^2* u(\tau)
 \rangle\bigr|\Bigr. \\
 +\Bigl. \bigl|\langle B\cdot\nabla B, e^{2\Delta(t-\tau)}\check{\varphi}^2
 * u(\tau)\rangle\bigr|\Bigr)~d\tau.
\end{multline}
\end{cor}

\noindent{\bf Proof.} Choose $E(t)=1$ and $\psi(\tau)$ as
\beqst
\psi(\tau)=\mathscr{F}^{-1}\bigl\{e^{-\abs{\xi}^2(t+\eta-\tau)}
\varphi(-\xi)\bigr\},\quad\eta>0
\eeqst
in Eq. (\ref{genener}). Then
$\psi(\tau)*u(\tau)=e^{\Delta(t+\eta-\tau)}\check{\varphi}*u(\tau)
=\int
e^{\Delta(t+\eta-\tau+\delta)}\check{\varphi}(x)u(\delta)~d\delta$,
and
\begin{align*}
\langle \psi'(\tau) & *u(\tau),\psi(\tau)*u(\tau)\rangle
-\norm{\nabla\psi(\tau)*u(\tau)}^2_2 \\
  & =-\bigl\langle \Delta\bigl(e^{\Delta(t+\eta-\tau)}\check{\varphi}\bigr)
    *u(\tau),e^{\Delta(t+\eta-\tau)}\check{\varphi}*u(\tau) \bigr\rangle
    -\norm{\nabla\bigl(e^{\Delta(t+\eta-\tau)}\check{\varphi}\bigr)*u(\tau)}^2_2
    \\
  & =0.
\end{align*}
Hence we have from (\ref{genener})
\begin{align*}
\norm{e^{\Delta\eta}\check{\varphi}*u(t)}_2^2 & \le
\norm{e^{\Delta(t+\eta-s)}\check{\varphi}*u(s)}^2_2 \\
 & \phantom{\le}~+2\int^t_s \Bigl( \abs{\langle u\cdot\nabla u,
    e^{2\Delta(t+\eta-\tau)}\check{\varphi}^2*u(\tau)\rangle}\Bigr.
    \\
    & \phantom{~+2\int^t_s va}+ \Bigl.\abs{\langle B\cdot\nabla B,
    e^{2\Delta(t+\eta-\tau)}\check{\varphi}^2*u(\tau)\rangle}\Bigr)~d\tau
\end{align*}
Let $\eta\to 0$ in the preceding equation to obtain
(\ref{mylab}).\hfill$\square$


\begin{cor}\label{corgen2}
Let $E(t)\in C^1(\R^+;\R)$ and $\tilde{\psi}\in
C^1(0,\infty;L^{\infty}\cap L^2)$. Then a weak solution of Equations
(\ref{mhd1})-(\ref{mhd4}) satisfies
\beq\label{eqcor}
\begin{split}
E(t)\bigl\|\tilde{\psi}(t)\Hat{u}(t)\bigr\|^2_2 & \le
E(s)\bigl\|\tilde{\psi}(s)\hat{u}(s)\bigr\|^2_2 + \int^t_s E'(\tau)
\bigl\|\tilde{\psi}(\tau)\hat{u}(\tau)\bigr\|^2_2~d\tau \\
  & + 2\int^t_s E(\tau)\Bigl( \langle \tilde{\psi}'(\tau)
   \hat{u}(\tau),\overline{\tilde{\psi}(\tau)\hat{u}(\tau)}\rangle -
   \bigl\|\xi\tilde{\psi}(\tau)\hat{u}(\tau)\bigr\|^2_2\Bigr)~d\tau \\
  & - 2\int^t_s E(\tau)\Bigl( \langle \mathscr{F}\{u\cdot
  \nabla u (\tau) \},\tilde{\psi}^2\hat{u}(\tau) \rangle\Bigr.\\
  &\phantom{abbadoa}+ \Bigl.\langle \mathscr{F}\{ B\cdot\nabla B (\tau) \},
   \tilde{\psi}^2\hat{u}(\tau) \rangle \Bigr)~d\tau,
\end{split}
\eeq
for almost all $s\ge 0$ and all $t\ge s$.
\end{cor}

\noindent{\bf Proof.} Apply Plancherel's theorem to (\ref{genener}).
\hfill$\square$

\section{Non-uniform decay of solutions to the MHD equations with no magnetic diffusion}

\setcounter{teo}{0}

In this section it is shown that if there were to exist strong
bounded solutions in 3D, they can not have compensatory
oscillations. We analyze separately the energy of the high and low
frequencies of the solutions. The main tool for the analysis of the
high frequency is Fourier Splitting, see \cite{maria85}.

We first establish Theorem \ref{teore1} of the introduction, which
we recall for completeness.

\begin{teo}\label{miteo1}
Let $n=3$. Suppose that there exists $(u,B)$ a strong bounded
solution to the MHD equations with $\delta=0$, magnetic field
bounded in $L^{\infty}$, and data $(u_0,B_0)\in (L^1\cap L^2\times
L^2\cap L^{\infty})$. Then
\beqst
\norm{u(t)}_2\to 0,\quad\text{and}\quad \norm{B(t)}_2\to M,
\eeqst
as $t\to\infty$, where $M$ is some positive constant.
\end{teo}

\noindent{\bf Proof.} Split the velocity of the solution into low
and high frequency parts
\beqst
\norm{u(t)}_2=\norm{\Hat{u}(t)}_2\le \norm{\varphi \Hat{u}}_2
+\norm{(1-\varphi)\Hat{u}}_2,
\eeqst
where $\varphi$ is a function in Fourier space to be chosen
appropriately, to emphasize the low and high frequency of $u$.

\noindent\textbf{\textsf{Low frequency Decay}.} Set
$\varphi(\xi)=e^{-\abs{\xi}^2t}$, using the result of Corollary
\ref{corgen} and Plancherel theorem,
\begin{align*}
\norm{\varphi \Hat{u}(t)}^2 & \le
\norm{e^{-\abs{\xi}^2(t-s)}\varphi\Hat{u}(s)}^2 +2\int^t_s
\abs{\langle u\cdot\nabla u, e^{2\Delta(t-\tau)}\check{\varphi}^2* u
 \rangle}~d\tau \\
 &\phantom{\le}~+2\int^t_s \abs{\langle B\cdot\nabla B, e^{2\Delta(t-\tau)}\check{\varphi}^2
 * u\rangle}~d\tau \\
 & \le \norm{e^{-\abs{\xi}^2(t-s)}\varphi\Hat{u}(s)}_2^2 + 2\int^t_s
 \abs{\langle \check{\varphi}^2* u\cdot\nabla u,
 e^{2\Delta(t-\tau)}u\rangle}~d\tau \\
 &\phantom{\le} ~+2\int^t_s \abs{\langle \check{\varphi}^2* B\cdot\nabla B,
 e^{2\Delta(t-\tau)}u\rangle}~d\tau
\end{align*}
Clearly the first term on the right hand side satisfies
\beq\label{ecuateo1}
\limsup\limits_{t\to\infty}\bigl\|e^{-\abs{\xi}^2(t-s)}\varphi\Hat{u}(s)
\bigr\|_2^2=0.
\eeq
To bound the third term note first that
\begin{align*}
\norm{\check{\varphi}^2* B\cdot\nabla B}^2_2 & =
\sum_j\int\abs{\check{\varphi}^2*(B\cdot\nabla)B^j}^2~dx \le
C\sum_j\left(\sum_i\norm{\partial_i\check{\varphi}^2* B^iB^j}_2
    \right)^2 \\ & \le
C\norm{B}_2^2\left(\sum_i\norm{\partial_i\check{\varphi}^2}_2
 \right)^2,
\end{align*}
where $C$ is a positive constant.

Thus
\beqst
\abs{\langle \check{\varphi}^2* B\cdot\nabla
B,e^{2\Delta(t-\tau)}u\rangle}\le
C\norm{B(\tau)}^2_2\norm{u(\tau)}_2\left(
\sum_i\norm{\partial_i\check{\varphi}^2}_2\right)
\eeqst
Since $\check{\varphi}^2=(4\pi
t)^{-3}e^{-\nicefrac{\abs{x}^2}{2t}}$, it follows that
\beq\label{neweq}
\norm{\partial_i\check{\varphi}^2}_2^2 \le
C\int_{\R^3}\left(t^{-4}\abs{x}e^{-\nicefrac{\abs{x}^2}{2t}}
\right)^2 dx \le C
t^{-11/2}\int_{\R^3}\dfrac{e^{-\nicefrac{\abs{x}^2}{2t}}}{t^{3/2}}dx
= C t^{-11/2}.
\eeq
Hence,
\beqst
\abs{\langle\check{\varphi}^2* B\cdot\nabla B,e^{2\Delta(t-\tau)}u
\rangle} \le C\norm{B(\tau)}_2\norm{u(\tau)}_2
t^{-\nicefrac{11}{2}}\le C t^{-\nicefrac{11}{2}},
\eeqst
since the $L^2$-norm of $u$ and $B$ are bounded by the initial data.

Similarly, $\abs{\langle\check{\varphi}^2* u\cdot\nabla
u,e^{2\Delta(t-\tau)}u \rangle}\le C t^{-\nicefrac{11}{2}}$. Hence
by (\ref{ecuateo1}), (\ref{neweq}), and the last two inequalities
\beqst
\norm{\varphi\check{u}(t)}^2_2\to 0 \quad\text{as}\quad t\to\infty.
\eeqst

\noindent\textbf{\textsf{High frequency Decay}.} We will show that
\beq\label{higheq}
\lim_{t\to\infty}\norm{(1-\varphi)u}_2 \le \eps_0
\eeq
for all $\eps_0>0$. The Fourier splitting method will be used. Let
$\chi(\eps)=\{\xi: \abs{\xi}\le G(\eps) \}$, a neighborhood of the
origin, were $G$ will be specified below. Set
$\tilde{\psi}=1-\varphi$, where $\varphi$ is given above. Then
$\tilde{\psi}'=\abs{\xi}^2\varphi$, and Corollary (\ref{corgen2})
yields
\beq\label{highfeq}
\begin{split}
E(t)\norm{(1-\varphi)\Hat{u}(t)}^2_2 & \le
E(s)\norm{(1-\varphi)\Hat{u}(s)}^2_2 + \int^t_s E'(\tau)
\int_{\chi(\eps)}\abs{(1-\varphi)\Hat{u}(\tau)}^2~d\xi\,d\tau \\
 & +\int^t_s E'(\tau)\int_{R^3\backslash\chi(\eps)}
 \abs{(1-\varphi)\Hat{u}(\tau)}^2~d\xi\,d\tau \\
 & -2\int^t_s E(\tau)\norm{\xi(1-\varphi)\Hat{u}(\tau)}^2_2~d\tau \\
  &+2\int^t_s E(\tau)\langle\abs{\xi}^2 \varphi(\tau)\Hat{u}(\tau),(1-\varphi(\tau))
   \Hat{u}(\tau)\rangle~d\tau \\
 &+2\int^t_s E(\tau) \abs{\langle \widehat{u\cdot\nabla u},
   (1-\varphi)^2\Hat{u}(\tau)\rangle}d\tau \\
 &+2\int^t_s E(\tau) \abs{\langle \widehat{B\cdot\nabla B},
   (1-\varphi)^2\Hat{u}(\tau)\rangle}d\tau.
\end{split}
\eeq
The terms in the second and third row are bounded by
\beqst
\int^t_s
\bigl(E'(\tau)-2E(\tau)G^2(\eps)\bigr)\int_{\R^3\backslash\chi(\eps)}
\abs{(1-\varphi)\Hat{u}(\tau)}^2~d\xi\,d\tau.
\eeqst
Choose $E(t)=e^{\eps t}$ and $G(\eps)=\sqrt{\eps/2}$ hence
$E'(t)-2E(t)G^2(\eps)=0$, thus the above integral vanishes.

Divide Equation (\ref{highfeq}) by $E(t)$,
\beq\label{miecua}
\begin{split}
& \norm{(1-\varphi)\Hat{u}(t)}^2_2 \le
\dfrac{E(s)}{E(t)}\norm{(1-\varphi)\Hat{u}(s)}^2_2
+\dfrac{1}{E(t)}\int^t_s E'(\tau)
\int_{\chi(\eps)}\abs{(1-\varphi)\Hat{u}(\tau)}^2~d\xi\,d\tau \\
 & \phantom{v}+ \dfrac{2}{E(t)}\int^t_s E(\tau)\langle\abs{\xi}^2 \varphi(\tau)\Hat{u}(\tau),
   (1-\varphi(\tau))\Hat{u}(\tau)\rangle~d\tau \\
 & \phantom{v}+ \dfrac{2}{E(t)}\int^t_s E(\tau) \abs{\langle \widehat{u\cdot\nabla u},
   (1-\varphi)^2\Hat{u}(\tau)\rangle}d\tau
  + \dfrac{2}{E(t)}\int^t_s E(\tau) \abs{\langle \widehat{B\cdot\nabla B},
   (1-\varphi)^2\Hat{u}(\tau)\rangle}d\tau\\
 & \phantom{v}= I(t)+II(t)+III(t)+IV(t)+V(t).
\end{split}
\eeq

\noindent\textsc{Estimate for $I(t)$: }
\beqst
I(t)=\dfrac{E(s)}{E(t)}\norm{(1-\varphi)\Hat{u}(s)}^2_2 \le
e^{\eps(s-t)}\norm{\Hat{u}(s)}^2_2 \le
\Bigl(\norm{u_0}^2_2+\norm{B_0}^2_2 \Bigr) e^{\eps(s-t)}\le C
e^{\eps(s-t)}.
\eeqst
Hence
\beqst
\lim_{t\to\infty}I(t)=\lim_{t\to\infty}\dfrac{E(s)}{E(t)}\norm{(1-\varphi)\Hat{u}(s)}^2_2
=0.
\eeqst

\noindent\textsc{Estimate for $II(t)$: } By Proposition \ref{decay}
\beqst
\int_{\chi(\eps)} \abs{\Hat{u}}^2~d\xi  \le C\int_{\chi(\eps)}
\biggl(1+\dfrac{1}{\abs{\xi}} \biggr)^2 d\xi \le C\int_{\chi(\eps)}
\biggl( 1 +\dfrac{1}{\abs{\xi}^2} \biggr) d\xi
 \le C\bigl(\eps^{3/2}+\eps^{1/2} \bigr).
\eeqst
Since $(1-\varphi)^2\le 1$,
\beqst
\begin{split}
II(t) &= \dfrac{1}{E(t)}\int^t_s E'(\tau)
\int_{\chi(\eps)}\abs{(1-\varphi)\Hat{u}(\tau)}^2~d\xi\,d\tau \\
& \le \dfrac{1}{E(t)}\int^t_s E'(\tau)\int_{\chi(\eps)}
\abs{\Hat{u}(\tau)}^2 d\xi\, d\tau \le C\eps^{1/2}=\eps_0.
\end{split}
\eeqst

\noindent\textsc{Estimate for $III(t)$: } Observe that
$0\le\varphi-\varphi^2\le 1$, and $E(\tau)\le E(t)$ for $\tau<t$,
hence
\beqst
\begin{split}
III(t) & =\dfrac{2}{E(t)}\int^t_s E(\tau)\left\langle\abs{\xi}^2
\varphi(\tau)\Hat{u}(\tau),(1-\varphi(\tau))
\Hat{u}(\tau)\right\rangle d\tau \\
  &  \le \dfrac{2}{E(t)}\int^t_s E(t)\int_{\R^3}\abs{\xi}^2\bigl(
  \varphi-\varphi^2\bigr)\abs{\Hat{u}(\tau)}^2 d\xi\,d\tau
   \le C\int^t_s\norm{\nabla u(\tau)}^2_2 d\tau,
\end{split}
\eeqst
Since $\int^{\infty}_0\norm{\nabla u}_2^2<\infty$ it follows that
\beqst
\lim_{t\to\infty} III(t) \le C\lim_{s\to\infty}\lim_{t\to\infty}
\int^t_s\norm{\nabla u(\tau)}^2_2 d\tau=0
\eeqst

\noindent\textsc{Estimate for $IV(t)$: } Set
$\zeta=\mathscr{F}^{-1}\bigl\{1-(1-\varphi)^2 \bigr\}$. This
function is essentially the heat kernel. Note that $\langle
u\cdot\nabla u,u\rangle=0$, hence $IV(t)$ can be estimated as
follows
\begin{align*}
IV(t) &= \dfrac{2}{E(t)}\int^t_s E(\tau) \abs{\langle u\cdot\nabla
u,\zeta* u\rangle}~d\tau
  \le C\int^t_s \norm{\zeta}_{6/5}\norm{u\cdot\nabla u}_{3/2}
 \norm{u}_2 ~d\tau \\
 & \le C\int^t_s \norm{\zeta}_{6/5}\norm{u}_6\norm{\nabla u}_2~d\tau
 \le C\int^t_s \norm{\zeta}_{6/5}\norm{\nabla u}^2_2~d\tau \\
 & \le C\int^t_s \dfrac{\norm{\nabla u}_2^2}{\tau^{1/4}} d\tau
 \le \dfrac{C}{s^{1/4}} \int^t_s\norm{\nabla u}^2_2~d\tau.
\end{align*}
Thus
\beqst
\lim_{t\to\infty} IV(t)\le \lim_{s\to\infty}\lim_{t\to\infty}
\dfrac{C}{s^{1/4}} \int^t_s\norm{\nabla u}^2_2~d\tau=0.
\eeqst

\noindent\textsc{Estimate for $V(t)$: }

\begin{align*}
V(t) & = \dfrac{2}{E(t)}\int^t_s E(\tau) \abs{\langle
\widehat{B\cdot\nabla B},(1-\varphi)^2\Hat{u}(\tau)\rangle}d\tau \le
\dfrac{C}{E(t)}\int^t_s E(\tau)\left(\sum_j\int\abs{\widehat{B\cdot
B}}\abs{\widehat{\nabla u^j}} d\xi \right) d\tau \\
 & \le \dfrac{C}{E(t)}\int^t_s E(\tau) \left(\int\abs{B\cdot B}^2
 ~d\xi\right)^{\nicefrac{1}{2}}
 \left(\int\abs{\nabla u(\tau)}^2~d\xi
 \right)^{\nicefrac{1}{2}} d\tau \\
 & \le \dfrac{C}{E(t)}\int^t_s E(\tau) \left(\int \abs{B(\tau)}^4~d\xi\right)^{\nicefrac{1}{2}} \norm{\nabla
 u(\tau)}_2 d\tau \\
 & \le \dfrac{C}{E(t)}\int^t_s E(\tau) \norm{B}_{\infty} \norm{B(\tau)}_2
   \norm{\nabla u(\tau)}_2 d\tau \\
 & \le \dfrac{C}{E(t)}\int^t_s E(\tau) \norm{\nabla u(\tau)}_2
 d\tau,
\end{align*}
Here is the only place where we need that the magnetic field is
bounded. Specifically we use that $B\in L^{\infty}(\R^3\times\R_+)$.
Since $B\in L^2$ we get the above bound. Recall that $E(t)=e^{\eps
t}$, hence
\begin{align*}
V(t) & \le \dfrac{C}{E(\tau)}
    \int^t_s E(t)\norm{\nabla u(\tau)}_2 d\tau
   \le \dfrac{C}{e^{\eps t}}\left( \dfrac{e^{2\eps t}}{2\eps}
  \right)^{\nicefrac{1}{2}} \left(\int^t_s \norm{\nabla u(\tau)}^2_2
  d\tau\right)^{\nicefrac{1}{2}} \\
  & = \dfrac{C}{\sqrt{2\eps}} \left(\int^t_s \norm{\nabla u(\tau)}^2_2
  d\tau\right)^{\nicefrac{1}{2}},
\end{align*}
which as before tends to zero as $t$ and $s$ goes to infinity.

Combining the estimates $I(t)-V(t)$ yields
 In summary, we have showed that
\beq\label{limvarphi}
\lim_{t\to\infty}\norm{(1-\varphi)\Hat{u}(t)}_2\le \eps_0.
\eeq
Since $\eps_0$ is arbitrary and positive, combining (\ref{ecuateo1})
and (\ref{limvarphi}) yields
\beqst
\lim_{t\to\infty}\norm{u(t)}_2 =0.
\eeqst
To obtain the limit of $\norm{B}_2$ proceed as follows. Set
\beqst
\phi(t)=\norm{u(t)}^2_2 + \norm{B(t)}^2_2.
\eeqst
Given that $\phi(t)\ge 0$ and is decreasing, there exists a constant
$M$ such that $\phi(t)\to M$ as $t\to\infty$. Since
$\norm{u(t)}_2\to 0$, it follows that
\beqst
\norm{B(t)}_2\to M\qquad \text{as}\quad t\to\infty.
\eeqst
This completes the proof. \hfill$\square$

\section{MHD equations with diffusion}

In this section it is shown that if the data is only in $L^2$ then
the solution decays without a rate. The ideas of the proof are
similar to Theorem \ref{miteo1} only that due to the added magnetic
diffusion we need less information on the data. The main result of
this section is that this decay is optimal. Specifically it is shown
that the decay can not be uniform.

The proof we give below is formal. To make it rigorous, it can be
applied to smooth approximations and then pass to the limit. The
approximations could be constructed by retarded mollification as was
done for the Navier-Stokes equations in \cite{ckn82,leray34}. This
construction if modified for the MHD equations will give suitable
approximations which can be used to make our arguments rigorous.
This arguments are standard and as such will be omitted. To see the
construction of these approximations in detail we refer the reader
to \cite{xin05}.

\subsection{Non-uniform decay}

\begin{teo}\label{miteo3}
Let $n=2,3$. If $(u,B)$ is a weak solution to the MHD equations with
$\delta>0$, and data $(u_0,B_0)\in \bigl(L^2(\R^n)\bigr)^2$, then
\beqst
\lim_{t\to\infty}\left(\norm{u(t)}^2_2+\norm{B(t)}^2_2\right)= 0.
\eeqst
\end{teo}

\noindent\textbf{Proof.} Without loss of generality suppose
$\delta=1$. The proof is based on similar arguments given in
\cite{ors97} for solutions to the Navier-Stokes equations with a
forcing term.

Let $\phi(\xi)=e^{-\abs{\xi}^2}$. As before split $u$ into low and
high frequency parts
\beqst
\begin{split}
\norm{u(t)}_2 & =\norm{\Hat{u}(t)}_2 \le \norm{\varphi\Hat{u}}_2 +
\norm{(1-\varphi)\Hat{u}}_2, \\
\norm{B(t)}_2 & =\norm{\Hat{B}(t)}_2 \le \norm{\varphi\Hat{B}}_2 +
\norm{(1-\varphi)\Hat{B}}_2.
\end{split}
\eeqst
\noindent{\bf Low frequency decay.} We need to use Corollary
(\ref{corgen}) and Plancherel's identity.
\beq\label{nonuno}
\begin{split}
\norm{\varphi \Hat{u}(t)}^2 &  \le
\norm{e^{-\abs{\xi}^2(t-s)}\varphi\Hat{u}(s)}_2^2 + 2\int^t_s
 \abs{\langle \check{\varphi}^2* u\cdot\nabla u,
 e^{2\Delta(t-\tau)}u\rangle}~d\tau \\
 &\phantom{\le} ~+2\int^t_s \abs{\langle \check{\varphi}^2* B\cdot\nabla B,
 e^{2\Delta(t-\tau)}u\rangle}~d\tau=I_1+I_2+I_3,
\end{split}
\eeq
and
\beq\label{nondos}
\begin{split}
\norm{\varphi \Hat{B}(t)}^2 &  \le
\norm{e^{-\abs{\xi}^2(t-s)}\varphi\Hat{B}(s)}_2^2 + 2\int^t_s
 \abs{\langle \check{\varphi}^2* B\cdot\nabla B,
 e^{2\Delta(t-\tau)}u\rangle}~d\tau \\
 &\phantom{\le} ~+2\int^t_s \abs{\langle \check{\varphi}^2* u\cdot\nabla B,
 e^{2\Delta(t-\tau)}u\rangle}~d\tau=J_1+J_2+J_3.
\end{split}
\eeq
It is immediate that the first terms $I_1$ and $J_1$ in
(\ref{nonuno}), (\ref{nondos}) tend to zero respectively as $t$ goes
to infinity. Hence it will be only necessary to show that the two
integrals on the right hand side of each of the above equations
tends to zero when $t$ goes to infinity. Since all integrals can be
estimated in a similar fashion we will only analyze the integrals
$I_2$, and $I_3$ corresponding to the velocity.

Since $\Check{\varphi}^2$ is a rapidly decreasing function, by the
Hasdorff-Young, H\"{o}lder, and Sobolev inequalities we have the
following.

When $n=2$
\beqst
\begin{split}
\abs{\langle \check{\varphi}^2 * u\cdot \nabla u,e^{2\Delta
(t-\tau)}u\rangle} & \le
 \abs{\langle
 u\cdot\nabla\check{\varphi}^2*e^{2\Delta(t-\tau)}u,u\rangle}
  \le \norm{u}^2_4 \norm{\check{\varphi}^2* e^{2\Delta(t-\tau)}\nabla
 u}_2 \\
& \le \norm{\check{\varphi}^2}_{\infty} \norm{u}^2_4\norm{\nabla
u}_2
  \le C\norm{u}_2 \norm{\nabla u}^2_2
\end{split}
\eeqst
and
\beqst
\begin{split}
\abs{\langle \check{\varphi}^2 * B\cdot \nabla B,e^{2\Delta
(t-\tau)}u\rangle} & \le
 \abs{\langle
 B\cdot\nabla\check{\varphi}^2*e^{2\Delta(t-\tau)}B,u\rangle}
  \le \norm{B}^2_4 \norm{\check{\varphi}^2* e^{2\Delta(t-\tau)}\nabla
 u}_2 \\
& \le \norm{\check{\varphi}^2}_{\infty} \norm{B}^2_4\norm{\nabla
u}_2
  \le C\norm{B}_2 \norm{\nabla B}_2 \norm{\nabla u}_2 \\
  & \le C\left(\norm{\nabla B}^2_2+\norm{\nabla u}^2_2 \right).
\end{split}
\eeqst

For $n=3$
\beqst
\begin{split}
\abs{\langle \check{\varphi}^2 * u\cdot \nabla u,e^{2\Delta
(t-\tau)}u\rangle} & \le \norm{\check{\varphi}^2 *u\cdot\nabla u}_2
\norm{u}_2
  \le C \norm{\check{\varphi}^2}_{6/5}\norm{u\cdot\nabla
 u}_{3/2}\norm{u}_2 \\
 & \le C \norm{u}_{6}\norm{\nabla u}_2
 \le C\norm{u}_2 \norm{\nabla u}_2^2.
\end{split}
\eeqst
In the same fashion
\beqst
\begin{split}
\abs{\langle \check{\varphi}^2 * B\cdot \nabla B,e^{2\Delta
(t-\tau)}u\rangle} & \le \norm{\check{\varphi}^2 *B\cdot\nabla B}_2
\norm{u}_2
  \le C \norm{\check{\varphi}^2}_{6/5}\norm{B\cdot\nabla
 B}_{3/2}\norm{u}_2 \\
 & \le C \norm{B}_{6}\norm{\nabla u}_2
 \le C\norm{u}_2 \norm{\nabla B}_2^2.
\end{split}
\eeqst

Hence integrating over $(s,t)$ yields
\beqst
I_2+I_3 \le C\int_s^t \norm{\nabla B}^2_2+\norm{\nabla u}^2_2~d\tau.
\eeqst
Thus
\beqst
\lim_{t\to\infty}I_2 +I_3 \le \lim_{s\to\infty}\lim_{t\to\infty}
C\int_s^t \norm{\nabla B}^2_2+\norm{\nabla u}^2_2~d\tau=0.
\eeqst
In the same manner it follows that the $\lim_{t\to\infty}J_2+J_3=0$.
Hence
\beq\label{unotres}
\lim_{t\to\infty}\norm{\varphi\Hat{u}}_2=0,\quad
\lim_{t\to\infty}\norm{\varphi\Hat{B}}_2=0.
\eeq

\noindent{\bf High frequency decay.} To estimate the high frequency
part we will use Fourier Splitting \cite{maria85}. We now use
Corollary (\ref{corgen2}) and an equivalent version of this
corollary for the magnetic field. Choose
$\tilde{\psi}=1-e^{-\abs{\xi}^2}=1-\varphi$ (note that in this case
$\tilde{\psi}$ is independent of time). Let $\chi(t)=\{\xi\in\R^n:
\abs{\xi}\le G(t) \}$, then
\beq\label{uno4}
\begin{split}
E(t) &\left[\norm{(1-\varphi)\Hat{u}(t)}^2_2 +
\norm{(1-\varphi)\Hat{B}(t)}^2_2 \right] \le
E(s)\left(\norm{(1-\varphi)\Hat{u}(s)}^2_2+
\norm{(1-\varphi)\Hat{B}(s)}^2_2\right) \\ &+ \int^t_s E'(\tau)
\int_{\chi(t)}\abs{(1-\varphi)\Hat{u}(\tau)}^2~d\xi\,d\tau +
\int^t_s E'(\tau)
\int_{\chi(t)}\abs{(1-\varphi)\Hat{B}(\tau)}^2~d\xi\,d\tau\\
 & +\int^t_s E'(\tau)\int_{R^3\backslash\chi(t)}
 \abs{(1-\varphi)\Hat{u}(\tau)}^2~d\xi\,d\tau
 +\int^t_s E'(\tau)\int_{R^3\backslash\chi(t)}
 \abs{(1-\varphi)\Hat{B}(\tau)}^2~d\xi\,d\tau\\
 & -2\int^t_s E(\tau)\norm{\xi(1-\varphi)\Hat{u}(\tau)}^2_2~d\tau
 -2\int^t_s E(\tau)\norm{\xi(1-\varphi)\Hat{B}(\tau)}^2_2~d\tau\\
 &-2\int^t_s E(\tau) \langle \widehat{u\cdot\nabla u},
   (1-\varphi)^2\Hat{u}(\tau)\rangle d\tau
 +2\int^t_s E(\tau) \langle \widehat{B\cdot\nabla B},
   (1-\varphi)^2\Hat{u}(\tau)\rangle d\tau \\
 & -2\int^t_s E(\tau) \langle \widehat{B\cdot\nabla u},
   (1-\varphi)^2\Hat{B}(\tau)\rangle d\tau
 +2\int^t_s E(\tau) \langle \widehat{u\cdot\nabla B},
   (1-\varphi)^2\Hat{B}(\tau)\rangle d\tau.
\end{split}
\eeq
Suppose
\beq\label{eyg}
 E(t)=(1+t)^{\alpha}, \quad\text{and}\quad
G^2(t)=\tfrac{\alpha}{2(1+t)},
\eeq
with $\alpha>3$. With this choice $E'(t)-2E(t)G^2(t)=0$, proceeding
as in Theorem \ref{miteo1} will yield
\beq\label{uno0}
\begin{split}
\int^t_s & E'(\tau)\int_{R^3\backslash\chi(t)}
 \abs{(1-\varphi)\Hat{u}(\tau)}^2~d\xi\,d\tau
-2\int^t_s E(\tau)\norm{\xi(1-\varphi)\Hat{u}(\tau)}^2_2~d\tau\\
& +\int^t_s E'(\tau)\int_{R^3\backslash\chi(t)}
 \abs{(1-\varphi)\Hat{B}(\tau)}^2~d\xi\,d\tau
 -2\int^t_s E(\tau)\norm{\xi(1-\varphi)\Hat{B}(\tau)}^2_2~d\tau \le
 0.
\end{split}
\eeq
By (\ref{uno0}), Equation (\ref{uno4}) can be reduced to
\begin{align}
E(t) &\left[\norm{(1-\varphi)\Hat{u}(t)}^2_2 +
\norm{(1-\varphi)\Hat{B}(t)}^2_2 \right] \le
E(s)\left(\norm{(1-\varphi)\Hat{u}(s)}^2_2+
\norm{(1-\varphi)\Hat{B}(s)}^2_2\right) \label{unoeq}\\
&+ \int^t_s E'(\tau)
\int_{\chi(t)}\abs{(1-\varphi)\Hat{u}(\tau)}^2~d\xi\,d\tau +
\int^t_s E'(\tau)
\int_{\chi(t)}\abs{(1-\varphi)\Hat{B}(\tau)}^2~d\xi\,d\tau\label{doseq}\\
&-2\int^t_s E(\tau) \langle \widehat{u\cdot\nabla u},
   (1-\varphi)^2\Hat{u}(\tau)\rangle d\tau
 +2\int^t_s E(\tau) \langle \widehat{B\cdot\nabla B},
   (1-\varphi)^2\Hat{u}(\tau)\rangle d\tau \label{treseq}\\
 & -2\int^t_s E(\tau) \langle \widehat{B\cdot\nabla u},
   (1-\varphi)^2\Hat{B}(\tau)\rangle d\tau
 +2\int^t_s E(\tau) \langle \widehat{u\cdot\nabla B},
   (1-\varphi)^2\Hat{B}(\tau)\rangle d\tau. \label{cuaeq}
\end{align}
We now will bound the terms in (\ref{doseq}). Observing that
$\abs{1-\varphi}\le \abs{\xi}^2$ if $\abs{\xi}<1$, we have
\beq\label{cge}
\begin{split}
\int_{\chi(\tau)}
\abs{(1-\varphi)}^2\left(\abs{\Hat{u}}^2+\abs{\Hat{B}}^2\right) d\xi
&\le C G(\tau)^4 \int_{\chi(\tau)}
\left(\abs{\Hat{u}}^2+\abs{\Hat{B}}^2\right) d\xi \\
 &\le C \left(
\norm{u_0}^2_2+\norm{B_0}^2_2\right) (1+\tau)^{-2}
\end{split}
\eeq

We now analyze (\ref{treseq}), (\ref{cuaeq}) together. For this note
first that $(1-\varphi)^2=1+\theta$, where
$\theta=-2\varphi+\varphi^2$, hence by the definition of $\varphi$,
the function $\theta$ is a rapidly decreasing function.

Since $\langle u\cdot\nabla u,u\rangle=\langle u\cdot\nabla B,B
\rangle=0$, and $\langle \widehat{B\cdot\nabla B},\Hat{u}\rangle-
\langle \widehat{u\cdot\nabla B},\Hat{B}\rangle=0$, it follows that
the four last terms of the right hand side of (\ref{uno4}) can be
expressed as
\beqst
\begin{split}
-2\int^t_s &E(\tau)  \langle \widehat{u\cdot\nabla u},
   \theta\Hat{u}(\tau)\rangle d\tau
 +2\int^t_s E(\tau) \langle \widehat{B\cdot\nabla B},
   \theta\Hat{u}(\tau)\rangle d\tau \\
 & -2\int^t_s E(\tau) \langle \widehat{B\cdot\nabla u},
   \theta\Hat{B}(\tau)\rangle d\tau
 +2\int^t_s E(\tau) \langle \widehat{u\cdot\nabla B},
   \theta\Hat{B}(\tau)\rangle d\tau
   = K_1+K_2+K_3+K_4.
\end{split}
\eeqst
The estimates of $K_i$'s are all very similar. Hence we only
estimate $K_1$ and state estimates are for $K_i$ for $i=2,3,4$.

For $n=2$
\beqst
\begin{split}
K_1 &= \int^t_s E(\tau)\abs{u\cdot\nabla u,\check{\theta}* u(\tau)}
d\tau
 = \int^t_s E(\tau)\abs{\langle u\cdot \check{\theta} *\nabla u(\tau),u(\tau)
 \rangle} d\tau \\
 & =\int^t_s E(\tau) \norm{\check{\theta}}_1 \norm{u}^2_4 \norm{\nabla u}_2
 d\tau
  \le C\int^t_s E(\tau)\norm{u}_2 \norm{\nabla u}^2_2 d\tau \\
 & \le C\left(\norm{u_0}^2_2
+\norm{B_0}^2_2\right)\int^t_s E(\tau)\norm{\nabla u}^2_2 d\tau,
\end{split}
\eeqst
Same type of computations yields that
\beqst
K_i \le C \int^t_s E(\tau)\left(\norm{\nabla u}^2_2+\norm{\nabla
B}^2_2 \right) d\tau,\quad i=2,3,4.
\eeqst

For $n=3$ we also only estimate $K_1$
\beqst
\begin{split}
\int^t_s E(\tau)\abs{\langle u\cdot \nabla u,\check{\theta}*
u(\tau)\rangle} d\tau &\le \int^t_s
E(\tau)\norm{\check{\theta}}_{6/5} \norm{u\cdot \nabla u}_{3}
\norm{u}_2 d\tau \\
 & \le C \int^t_s E(\tau) \norm{u}_{6}\norm{\nabla u}_2 d\tau
  \le C\int^t_s E(\tau) \norm{\nabla u}^2_2 d\tau.
\end{split}
\eeqst
The same estimates yield
\beqst
K_i \le C \int^t_s E(\tau)\left(\norm{\nabla u}^2_2+\norm{\nabla
B}^2_2 \right) d\tau,\quad i=2,3,4.
\eeqst
Combining the estimates (\ref{unoeq})--(\ref{cuaeq}), and the
estimates for the $K_i$'s yields after division by $E(t)$
\beqst
\begin{split}
\norm{(1-\varphi)\Hat{u}(t)}^2_2  &+
\norm{(1-\varphi)\Hat{B}(t)}^2_2  \le \dfrac{E(s)}{E(t)}
\norm{(1-\varphi)\Hat{u}(s)}_2^2 \\ &+ \dfrac{1}{E(t)}\int^t_s
E'(\tau) \int_{\chi(\tau)} \abs{1-\varphi}^2
\left(\abs{\Hat{u}}^2 +\abs{\Hat{B}}^2 \right) d\xi d\tau \\
 & + \dfrac{C}{E(t)}\int^t_s E(\tau) \left(\norm{\nabla u(\tau)}^2_2 +
 \norm{\nabla B(\tau)}^2_2 \right) d\tau.
\end{split}
\eeqst
Since $1-\varphi\le 1$, combining the last equation with
(\ref{cge}), recalling the definition of $E(t)$ and $G(t)$ in
(\ref{eyg}), and since $\alpha>3$ we have
\beqst
\begin{split}
\lim_{t\to\infty} &\left( \norm{(1-\varphi)\Hat{u}(t)}^2_2 +
\norm{(1-\varphi)\Hat{B}(t)}^2_2\right) \le \lim_{t\to\infty}
\left(\dfrac{1+s}{1+t} \right)^{\alpha} \left(\norm{u_0}^2_2
+\norm{B_0}^2_2\right) \\
& + C\left(\norm{u_0}^2_2 +\norm{B_0}^2_2\right)\lim_{t\to\infty}
 \left( \dfrac{1}{(1+t)^{\alpha}}\int^t_s (1+t)^{\alpha-3}\right)
 d\tau \\
 &+ \lim_{t\to\infty}\dfrac{C}{E(t)} \int^t_s E(\tau) \left(\norm{\nabla u}^2_2 +
   \norm{\nabla B}^2_2 \right) d\tau \\
  & = C\int^{\infty}_s \left(\norm{\nabla u}^2_2 +
   \norm{\nabla B}^2_2 \right) d\tau
\end{split}
\eeqst
Letting $s\to\infty$ on the right hand side yields
\beqst
\lim_{t\to\infty} \left( \norm{(1-\varphi)\Hat{u}(t)}^2_2 +
\norm{(1-\varphi)\Hat{B}(t)}^2_2\right) =0
\eeqst
Combining (\ref{unotres}) with the last limit gives
\beqst
\lim_{t\to\infty}\left(\norm{u}^2_2+\norm{B}^2_2 \right)=0.
\eeqst
As stated in the beginning to make this proof rigorous, the formal
proof has to be applied to the approximating solutions described at
the beginning of the section, and then pass to the limit. This
procedure is standard and as such is omitted. This completes the
proof of the theorem.\hfill$\square$

\subsection{Lack of uniformity}

In this subsection it is shown that solutions with data $(u_0,B_0)$
in $\bigl(L^2(\R^n)\bigr)^2$ cannot be expected to decay at a
uniform rate, in the sense that for each sphere in $L^2(\R^n)$ of
radius $\beta$, there is a point on the sphere so that the solution
with such data will decay arbitrarily slow. In other words, given a
time $T>0$, and $\eps>0$, there exists data $u_0$ with
$\norm{u_0}_2=\beta$ so that
\beq\label{ratio1}
\dfrac{\norm{u(T)}_2}{\norm{u_0}_2}\ge 1-\eps.
\eeq
Similarly for $B$ given a time $T>0$, and $\eps>0$, there exists
data $B_0$ with $\norm{B_0}_2=\gamma$ such that
\beq\label{btrece}
\dfrac{\norm{B(T)}_2}{\norm{B_0}_2}\ge 1-\eps.
\eeq

\begin{teo}
There exist no functions $G(t,\beta)$ and $H(t,\gamma)$ with the
following two properties. If $(u,B)$ is a solution to equations
(\ref{mhd1})-(\ref{mhd4}) with $\delta>0$ and data $(u_0,B_0)\in
(L^2(\R^n))^2$, $n=2,3$, then
\begin{enumerate}[i{)}]
\item $\norm{u(t)}_2\le G(t,\norm{u_0}_2)$, and $\norm{B(t)}_2\le
H(t,\norm{B_0}_2)$,
\item $\lim_{t\to 0}G(t,\beta)=0$, for all $\beta >0$, and $\lim_{t\to 0}
H(t,\gamma)=0$, for all $\gamma>0$.
\end{enumerate}
\end{teo}

\noindent\textbf{Proof.} As pointed in \cite{maria86} this lack of
uniformity is already present at the level of the heat equation.

The proof of Proposition 2.1 in \cite{maria86} has a gap that will
be taken care of in our present work.

Notice that it suffices to show that (\ref{ratio1}) and
(\ref{btrece}) hold. The plan is to choose as initial data a family
$\{u_0^{\alpha},B_0^{\alpha} \}$ which satisfy
\beq\label{igualnorm}
\norm{u_0^{\alpha}}_2=\norm{u_0}_2\quad\text{and}
\quad\norm{B_0^{\alpha}}_2 =\norm{B_0}_2.
\eeq
Write the solutions in integral form
\beqst
\begin{split}
u^{\alpha}(x,t) & = K*u_0^{\alpha}-\int^t_0
K(x-y,t-s)*(u^{\alpha}\nabla
 u^{\alpha}-B^{\alpha}\nabla B^{\alpha}+\nabla p^{\alpha})ds, \\
B^{\alpha}(x,t) & = K*B_0^{\alpha}-\int^t_0
K(x-y,t-s)*(u^{\alpha}\nabla B^{\alpha}-B^{\alpha}\nabla u^{\alpha})
ds,
\end{split}
\eeqst
then
\beq\label{partab}
\begin{split}
\norm{u^{\alpha}(x,t)}_2 & \ge \norm{K*u_0^{\alpha}}_2-\int^t_0
\norm{K(x-y,t-s)*(u^{\alpha}\nabla
 u^{\alpha}-B^{\alpha}\nabla B^{\alpha}+\nabla p^{\alpha})}_2 ds, \\
\norm{B^{\alpha}(x,t)}_2 & \ge \norm{K*B_0^{\alpha}}_2-\int^t_0
\norm{K(x-y,t-s)*(u^{\alpha}\nabla B^{\alpha}-B^{\alpha}\nabla
u^{\alpha})}_2 ds.
\end{split}
\eeq
The choice of $(u_0^{\alpha},B_0^{\alpha})$ will be the following
\beqst
u_0^{\alpha}(x)=\alpha^{n/2} u_0(x),\quad
B_0^{\alpha}(x)=\alpha^{n/2} B_0(x),\quad\text{for }n=2,3.
\eeqst
Then it is easy to see that $u_0^{\alpha}, B_0^{\alpha}$ are
invariant under $\alpha$ in $L^2$ (i.e., Equation (\ref{igualnorm})
is satisfied). Hence it is only necessary to show
\beqst
\dfrac{\norm{K*u^{\alpha}_0}}{\norm{u_0}}\ge
1-\eps,\quad\text{and}\quad
\dfrac{\norm{K*B^{\alpha}_0}}{\norm{B_0}}\ge,
\eeqst
and
\beqst
\begin{split}
\int^t_0 \norm{K(x-y,t-s)*(u^{\alpha}\nabla
 u^{\alpha}-B^{\alpha}\nabla B^{\alpha}+\nabla p^{\alpha})}_2 ds<
 \tilde{\eps}, \\
\int^t_0 \norm{K(x-y,t-s)*(u^{\alpha}\nabla
B^{\alpha}-B^{\alpha}\nabla u^{\alpha})}_2 ds< \tilde{\eps},
\end{split}
\eeqst
for $\tilde{\eps}$ sufficiently small.

We also note that the data $(u_0^{\alpha},B_0^{\alpha})$ will yield
for the linear part a self-similar solution, that is
\beq\label{selfsim}
\begin{split}
K*u_0^{\alpha}(x,t) &= \alpha^{n/2}\Tilde{u}(\alpha x, \alpha^2 t), \\
K*B_0^{\alpha}(x,t) &= \alpha^{n/2}\Tilde{B}(\alpha x, \alpha^2 t),
\end{split}
\eeq
Equalities (\ref{selfsim}) follow by uniqueness of the heat equation
and since $\alpha^{n/2}\Tilde{u}(\alpha x, \alpha^2 t)$ and
$K*u_0^{\alpha}(x,t)$ are solutions to the heat equation
$u^{\alpha}_0$. In the same way $\alpha^{n/2}\Tilde{B}(\alpha x,
\alpha^2 t)$ and $K*B_0^{\alpha}(x,t)$ are solutions of the heat
equation with data $B_0^{\alpha}$.
We only will show the proof for the energy of the velocity since the
proof for the energy of the magnetic field is almost identical.

We show first that Equation (\ref{ratio1}) holds for the data
$u_0^{\alpha}$.
\begin{align*}
\int_{\R^n}\abs{\Tilde{u}^{\alpha}}^2 dx &=\alpha^n\int_{\R^n}
\abs{\Tilde{u}(\alpha x, \alpha^2 t)}^2
dx=\int_{\R^n}\abs{\Tilde{u}(y,\alpha^2 t)}^2 dy\\
 & =\int_{\R^n} \abs{\Hat{\Tilde{u}}(\xi,\alpha^2 t)}^2 d\xi =
 \int_{\R^n} e^{-2\abs{\xi}^2\alpha^2 t}\abs{\Hat{u}_0(\xi)}^2 d\xi
\end{align*}
Next, by Lebesgue Dominated Convergence theorem it follows that for
each fixed $t$,
\beqst
\lim_{\alpha\to 0}\int_{\R^n} e^{-2\abs{\xi}^2\alpha^2
t}\abs{\Hat{u}_0(\xi)}^2 d\xi=\int_{\R^n}\abs{\Hat{u}_0}^2 d\xi.
\eeqst
Hence
\beq\label{ratio}
\lim_{\alpha\to
0}\dfrac{\norm{\Tilde{u}^{\alpha}(\cdot,t)}^2_2}{\norm{u_0(\cdot)}^2_2}=1,
\eeq
Now it is necessary to show
\beq\label{foru}
\lim_{\alpha\to 0}\int^t_0 \norm{K(x-y,t-s)*\bigl(u^{\alpha}\nabla
 u^{\alpha}-B^{\alpha}\nabla B^{\alpha}+\nabla p^{\alpha}\bigr)}_2 ds=0,
\eeq
Hence we analyze
\beq\label{unomas}
\begin{split}
\int^T_0 \left\| K(x-y,t-s)\right. & \left. * \,
\bigl(u^{\alpha}\nabla
 u^{\alpha}-B^{\alpha}\nabla B^{\alpha}+\nabla p^{\alpha}\bigr)\right\|_2 ds \\
 & \le \int^T_0 \norm{\nabla K(t-s)}_2\biggl(\norm{u^{\alpha} u^{\alpha}}_2
   +\norm{B^{\alpha} B^{\alpha}}_2 +\norm{p^{\alpha}}_2
   \biggr)ds \\
 & \le C\int^T_0\dfrac{1}{(t-s)^{1/2}}\biggl(\norm{u^{\alpha} u^{\alpha}}_2
   +\norm{B^{\alpha} B^{\alpha}}_2 \biggr)ds.
\end{split}
\eeq
Here we used that
\beqst
\norm{p^{\alpha}}_2 \le C\biggl(\norm{u^{\alpha} u^{\alpha}}_2
+\norm{B^{\alpha} B^{\alpha}}_2 \biggr).
\eeqst
We can suppose then that we have chosen $u_0$ and $B_0$ to be in
$H^1$. Now observe that for $n=3$
\beqst
\norm{u^{\alpha} u^{\alpha}}_2 \le C\norm{\nabla
u^{\alpha}}^{3/2}_2,\quad \norm{B^{\alpha} B^{\alpha}} \le
C\norm{\nabla B^{\alpha}}^{3/2}_2
\eeqst
and for $n=2$
\beqst
\norm{u^{\alpha} u^{\alpha}}_2 \le C\norm{\nabla u^{\alpha}}_2,\quad
\norm{B^{\alpha} B^{\alpha}}_2 \le C\norm{\nabla B^{\alpha}}_2.
\eeqst
Also
\begin{align*}
\norm{\nabla u}_2 +\norm{\nabla B}_2 & \le C\bigl( \norm{\nabla
u}^2_2 +\norm{\nabla B}^2_2\bigr)^{1/2} \quad\text{and}\\
\norm{\nabla u}^{3/2}_2 +\norm{\nabla B}^{3/2}_2 & \le C \bigl(
\norm{\nabla u}^2_2 +\norm{\nabla B}^2_2\bigr)^{3/4}.
\end{align*}
Thus to bound the right hand side of (\ref{unomas}) we need to
estimate
\beqst
\varphi(t)=\norm{\nabla u^{\alpha}}^2_2+\norm{\nabla
B^{\alpha}}^2_2.
\eeqst
Note first that if we choose $(u_0,B_0)\in H_1\times H_1$ will yield
that $\norm{\nabla u^{\alpha_0}}_2$, $\norm{\nabla B^{\alpha_0}}_2$
are arbitrarily small if $\alpha\ll 1$. This follows since
\begin{align*}
\varphi(0) &=\varphi_0 =\int \abs{\nabla u^{\alpha}_0}^2 dx +\int
\abs{\nabla B^{\alpha}_0}^2 dx \\
 & = \int \alpha^4 \abs{\nabla u_0(\alpha x)}^2 dx +\int \alpha^4\abs{\nabla B_0(\alpha x)}^2
 dx \\
 & =\alpha^2 \left(\int \abs{\nabla u_0(x)}^2 dx +\int\abs{\nabla B_0(x)}^2 dx\right) \\
 & = C \alpha^2.
\end{align*}
In order to estimate $\norm{\nabla u^{\alpha}}_2$ and $\norm{\nabla
B^{\alpha}}_2$ we use Prodi's inequality. We consider two cases:
\begin{description}
\item[Case 1] ($n=2$). By Prodi,
\beqst
\dfrac{d\varphi}{dt} \le C \varphi^2 \quad \Rightarrow\quad
\dfrac{d\varphi}{\varphi}\le C \varphi dt \quad \Rightarrow \quad
\ln\left(\dfrac{\varphi(t)}{\varphi^{\alpha}_0} \right)\le
C\int^t_0\varphi(s)~ds,
\eeqst
which implies
\begin{align*}
\varphi(t) & \le \varphi^{\alpha}_0 e^{C\int^t_0\varphi(s)~ds}\le
\varphi^{\alpha}_0
e^{C(\norm{u_0}^2_2+\norm{B_0}^2_2)}\le \varphi^{\alpha}_0 \\
  & \le C\alpha^2.
\end{align*}
Hence,
\beq
\int^t_0 \norm{K(x-y,t-s)*\bigl(u^{\alpha}\nabla
 u^{\alpha}-B^{\alpha}\nabla B^{\alpha}+\nabla p^{\alpha}\bigr)}_2 ds\le CT^{1/2}
 \alpha^2,
\eeq
and thus (\ref{foru}) follows when $n=2$.

By $\ref{foru}$ and (\ref{ratio}) the theorem follows for the
velocity in two dimensions.
\item[Case 2] ($n=3$). By Prodi,
\beqst
\dfrac{d\varphi}{dt} \le C\varphi^3\quad \Rightarrow \quad
\dfrac{d\varphi}{\varphi^2} \le C \varphi dt.
\eeqst
Integrating,
\begin{align*}
\dfrac{1}{\varphi_0}-\dfrac{1}{\varphi(t)} & \le
C\int^t_0\varphi(s)~ds \\
  & \le C \bigl(\norm{u_0}^2_2 +\norm{B_0}^2_2 \bigr)=C.
\end{align*}
Solving for $\varphi(t)$ we get
\beqst
\varphi(t)\le\dfrac{\varphi_0}{1-C_0\varphi_0}=\dfrac{C\alpha^2}{1-\Tilde{C}
\alpha^2} \le 2 C\alpha^2,
\eeqst
where we have chosen $\alpha$ very small, say
$\Tilde{C}\alpha^2\le\frac{1}{2}$, to make the last inequality true.
Hence, for the case $n=3$, the expression (\ref{foru}) is negligible
too. Combining (\ref{foru}) and (\ref{ratio}) yields the conclusion
of the theorem for the velocity in three dimensions.

The estimate for the magnetic field  in 2 and 3 dimensions follows
in an analogous fashion. This completes the proof of the theorem.
\end{description}
\hfill$\square$

\section{Kato's estimates}

In this section we show that
 by some simple modification
Kato's pioneering work on $L^p$ decay for Navier-Stokes equations
\cite{kato84} holds for the MHD equations with magnetic diffusion.
The main difference is that his approximating solutions will be
replaced by the corresponding ones from MHD. Thus, rewrite the MHD
equations as follows
\beqst
\partial_t v+Av+F(v)=0,
\eeqst
where $A=(A_1,A_2)=-P(\Delta,\delta\Delta)$, $v=(u,B)$,
$F(v)=F(v,v)=(F_1(v,v),F_2(v,v))$, and
\begin{align*}
F_1(u,w,B,D) &=P(u\nabla w)-P(B\nabla D) \\
F_2(u,w,B,D) &=P(u\nabla B)-P(B\nabla u).
\end{align*}
Here $P$ is the orthogonal projection of $L^2$ onto the subspace
$PL^2$, which denotes the collection of divergence-free elements of
$L^2$.

Following Kato
\begin{align*}
u_{n+1} &= u_0(t)+G_1(u_n,B_n), \\
B_{n+1} &= B_0(t)+G_2(u_n,B_n),
\end{align*}
where
\beqst
G_i(u,B)=-\int^t_0 e^{-(t-s)A_i}F_i\bigl(u(s),B(s)\bigr)\,ds
\eeqst
and
\beqst
(u_0(\cdot,t),B_0(\cdot,t))=\left(e^{-t A_1}u_0(x),e^{-t A_2}B_0(x)
\right).
\eeqst
Using these expressions appropriately in theorems 1,2, and 3 in
Kato's paper \cite{kato84} will yield the same results for the MHD
equations.

We want to show how these results can be used to extend the decay
results for the MHD equations in two dimensions when combined with
our $L^2$ results. We recall, for easy reference, Kato's first two
theorems. In this case $u$ stands for the solution to the
Navier-Stokes equations
\beqst
u_t + u\cdot \nabla u +\nabla p = \Delta u,\quad \diver u =0.
\eeqst

\begin{teononum}[Kato 1]
Let the initial data $a\in PL^m$. Then there is $T>0$ and a unique
solution $u$ such that
\begin{align}
 t^{(1-m/q)/2} u\in BC([0,T);PL^q) & & \text{for }m\le q \le \infty,
 \tag{1.1}\\
 t^{1-\nicefrac{m}{2q}}\nabla u\in BC([0,T);PL^q) & & \text{for }
 m\le q < \infty, \tag{1.1'}
\end{align}
both with values zero at $t=0$ except for $q=m$ in $(1.1)$, in which
$u(0)=a$. Moreover, $u$ has the additional property
\beq
u\in L^r((0,T_1);PL^q)\quad\text{with }1/r=(1-m/q)/2,\quad
m<q<m^2/(m-2), \tag{1.2}
\eeq
with some $0<T_1\le T$.
\end{teononum}

\begin{teononum}[Kato 2]
There is $\lambda>0$ such that if $\norm{a}_m\le\lambda$, then the
solution $u$ in Theorem (1) is global, i.e. we may take
$T=T_1=\infty$. In particular, $\norm{u(t)}_q$ decays like
$t^{-(1-m/q)/2}$ as $t\to\infty$, including $q=\infty$, and
$\norm{\nabla u(t)}_q$ decays like $t^{-(1-\nicefrac{m}{2q})}$,
including $q=m$.
\end{teononum}

As stated before usingg $(u_{n+1},B_{n+1})$ as defined above and
following Kato's proof with straightforward modifications yields

\begin{teo}\label{teore4}
Let $n=2,3$. Suppose $(u_0,B_0)\in\left(PL^p\cap PL^n(\R^n)
\right)^2$, where $1<p<n$. There exists $\lambda_1>0$ such that if
$\norm{u_0}_n\le\lambda_1$ and $\norm{B_0}_n\le\lambda_1$, then the
solution to the MHD equations with $\delta>0$ is global and for any
finite $q\ge p$
\beq\label{ecuateo4}
t^{(n/p-n/q)/2} (u,B)\quad\text{and}\quad t^{(n/p-n/q+1)/2} (\nabla
u, \nabla B) \in BC\bigl([1,\infty];PL^q \bigr)^2.
\eeq
\end{teo}

Combining the results of Theorems \ref{teore1} and \ref{teore4} in
the two dimensional case yields the following improved decay for the
solutions to the MHD equations with $\delta>0$.

\begin{cor}
There is $\lambda >0$ such that for $\norm{u_0}_2\le \lambda$ the
global solution of the equation (\ref{ecuateo4}) for $q\ge m,$ and
for $2\le r\le q$
\beqst
\lim_{t\to\infty} t^{\frac{r-2}{2r}}\norm{(u,B)}_r=0
\eeqst
\end{cor}

\noindent{\bf Proof.} It follows by interpolating $L^r$ between
$L^2$ and $L^q$ and using the decay rates of the solutions
corresponding to those Sobolev spaces.\hfill$\square$

\section{Appendix}

\begin{prop}\label{proofdecay}
 Let $(u,B)$ be a solution to the MHD equations (\ref{mhd1})-(\ref{mhd4}).
 Assume the initial data $u_0, B_0$ is in $L^1(\R^3)\cap L^2(\R^3)$.
 Then
 \beqst
\abs{\Hat{u}(t)}\le C\biggl(1+\dfrac{1}{\abs{\xi}} \biggr),
 \eeqst
 where $C$ is a constant.
\end{prop}

\noindent{\bf Proof.} We start by taking the Fourier transform of
Equation (\ref{mhd1})
\beqst
\Hat{u}_t + \widehat{u\cdot\nabla u}-\widehat{B\cdot\nabla B}
+\widehat{\nabla p} =-\abs{\xi}^2\Hat{u}.
\eeqst
Let us define
\beqst 
H(\xi,t) =\widehat{u\cdot\nabla u} -\widehat{B\cdot\nabla B}+
\widehat{\nabla p}.
\eeqst
Then $\Hat{u}_t+\abs{\xi}^2\Hat{u}=-H(\xi,t)$, and this equation can
be integrated using the method of integrating factor to get
\beqst
\Hat{u}(t)=\Hat{u}(0) e^{-\abs{\xi}^2t} -\int^t_0
H(\xi,s)e^{-\abs{\xi}^2(t-s)}~ds.
\eeqst
Then,
\beqst 
\abs{\Hat{u}(t)}\le\abs{\Hat{u}_0} +\int^t_0
\abs{H(\xi,s)}e^{-\abs{\xi}^2(t-s)}~ds.
\eeqst
To bound $\abs{H(\xi,s)}$, let us first bound $|\widehat{\nabla
p}|$. For this, let us take the divergence operator in Equation
(\ref{mhd1}) 
which yields
\beqst
\Delta p=\sum_{k,j}\dfrac{\partial^2}{\partial x^k\partial
x^j}(B^jB^k-u^ju^k).
\eeqst
It follows that
\begin{align*}
|\widehat{\nabla p} | & =\abs{\xi}\abs{\Hat{p}} \le \sum_{k,j}
\dfrac{\abs{\xi^k\xi^j}}{\abs{\xi}}\bigl(|\widehat{B^jB^k} | +
|\widehat{u^ju^k} | \bigr) \notag\\
  & \le \abs{\xi}\sum_{k,j} \bigl(|\widehat{B^jB^k} | +
|\widehat{u^ju^k} | \bigr)\label{u4}
\end{align*}
Hence
\begin{align*}
\abs{H(\xi,s)} &\le C\abs{\xi}\sum_{j,k}\bigl(|\widehat{u^ju^k}|
+|\widehat{B^jB^k}|
\bigr) \\
& \le C\abs{\xi}\bigl(\norm{u_0}^2_2 +\norm{B_0}^2_2\bigr) \le
C\abs{\xi}.
\end{align*}
Since $|\Hat{u}_0|\le\norm{u_0}_1=C$, it follows that
\begin{align*}
|\Hat{u}| &\le |\Hat{u}_0| +C\abs{\xi}\int^t_0 e^{-\abs{\xi}^2(t-s)
}ds\le C+ \dfrac{C}{\abs{\xi}}\Bigl(1-e^{-\abs{\xi}^2t} \Bigr) \\
 & \le C + \dfrac{C}{\abs{\xi}}= C\Bigl(1+\dfrac{1}{\abs{\xi}}
 \Bigr),
\end{align*}
which finishes the proof. \hfill$\square$


\bibliographystyle{alpha}
\bibliography{mibiblio}

\begin{thebibliography}{CKN82}

\bibitem[Cha81]{cha}
S.~Chandrasekhar.
\newblock {\em Hydrodynamic and Hydromagnetic Stability}.
\newblock Dover Publications, New York, 1981.

\bibitem[CKN82]{ckn82}
L.~Caffarelli, R.~Kohn, and L.~Nirenberg.
\newblock Partial regularity of suitable weak solutions of the
  {N}avier-{S}tokes equations.
\newblock {\em Comm. Pure Appl. Math.}, 35(6):771--831, 1982.

\bibitem[Cow76]{cow}
T.~Cowling.
\newblock {\em Magnetohydrodynamics}.
\newblock Monographs on astronomical subjects. Hilger, 2 edition, 1976.

\bibitem[HX05a]{xin052}
Cheng He and Zhouping Xin.
\newblock On the regularity of weak solutions to the magnetohydrodynamic
  equations.
\newblock {\em J. Differential Equations}, 213(2):235--254, 2005.

\bibitem[HX05b]{xin05}
Cheng He and Zhouping Xin.
\newblock Partial regularity of suitable weak solutions to the incompressible
  magnetohydrodynamic equations.
\newblock {\em J. Funct. Anal.}, 227(1):113--152, 2005.

\bibitem[Kat84]{kato84}
Tosio Kato.
\newblock Strong {$L\sp{p}$}-solutions of the {N}avier-{S}tokes equation in
  {${\R}\sp{m}$}, with applications to weak solutions.
\newblock {\em Math. Z.}, 187(4):471--480, 1984.

\bibitem[Kim02]{kim02}
Sangjeong Kim.
\newblock Gevrey class regularity of the magnetohydrodynamics equations.
\newblock {\em ANZIAM J.}, 43(3):397--408, 2002.

\bibitem[Koz87]{kozo87}
Hideo Kozono.
\newblock On the energy decay of a weak solution of the {MHD} equations in a
  three-dimensional exterior domain.
\newblock {\em Hokkaido Math. J.}, 16(2):151--166, 1987.

\bibitem[Ler34]{leray34}
J.~Leray.
\newblock Sur le mouvement d'un liquide visquex emplissant l'espace.
\newblock {\em Acta Mathematica}, 63:193--248, 1934.

\bibitem[LLP84]{landau}
L.D. Landau, E.M. Lifshitz, and L.P. Pitaevskii.
\newblock {\em Electrodynamics of Continuous Media}, volume~8 of {\em Course of
  theoretical physics}.
\newblock Butterworth-Heinemann, 2 edition, 1984.

\bibitem[MS89]{india89}
Satish~D. Mohgaonkar and R.~V. Saraykar.
\newblock {$L\sp 2$}-decay for solutions of the magnetohydrodynamic equations.
\newblock {\em J. Math. Phys. Sci.}, 23(1):35--55, 1989.

\bibitem[ORS97]{ors97}
Takayoshi Ogawa, Shubha~V. Rajopadhye, and Maria~E. Schonbek.
\newblock Energy decay for a weak solution of the {N}avier-{S}tokes equation
  with slowly varying external forces.
\newblock {\em J. Funct. Anal.}, 144(2):325--358, 1997.

\bibitem[OT00]{oliver00}
Marcel Oliver and Edriss~S. Titi.
\newblock Remark on the rate of decay of higher order derivatives for solutions
  to the {N}avier-{S}tokes equations in {${\bf R}\sp n$}.
\newblock {\em J. Funct. Anal.}, 172(1):1--18, 2000.

\bibitem[Sch85]{maria85}
Maria~Elena Schonbek.
\newblock {$L^2$} decay for weak solutions of the {N}avier-{S}tokes equations.
\newblock {\em Arch. Rational Mech. Anal.}, 88(3):209--222, 1985.

\bibitem[Sch86]{maria86}
Maria~E. Schonbek.
\newblock Large time behaviour of solutions to the {N}avier-{S}tokes equations.
\newblock {\em Comm. Partial Differential Equations}, 11(7):733--763, 1986.

\bibitem[SSS96]{maria96}
M.~E. Schonbek, T.~P. Schonbek, and Endre S{\"u}li.
\newblock Large-time behaviour of solutions to the magnetohydrodynamics
  equations.
\newblock {\em Math. Ann.}, 304(4):717--756, 1996.

\bibitem[Wu02]{wu02}
J.~Wu.
\newblock Bounds and new approaches for the 3{D} {MHD} equations.
\newblock {\em J. Nonlinear Sci.}, 12(4):395--413, 2002.

\end{thebibliography}

\end{document}